\documentclass[review]{elsarticle}

\usepackage{lineno,hyperref}
\modulolinenumbers[5]

\journal{Journal of Renewable Energy}

\usepackage{amsmath}
\usepackage{amssymb}
 \usepackage{amsthm}
 \usepackage{amsfonts}
\usepackage{multirow}
\usepackage{eqnarray}
\usepackage{verbatim}
\usepackage{soul}
\usepackage{diagbox} 
\usepackage{eqnarray}
\newcommand{\rP}{\mathrm{P\space }} 
\newcommand{\rE}{\mathrm{E\space }} 
\usepackage[ruled,vlined]{algorithm2e}
\usepackage{algorithmic}

\hypersetup{pdfauthor={Name}}
\begin{document}
\nolinenumbers
\begin{frontmatter}

\title{Optimal preventive maintenance scheduling for wind turbines under condition monitoring}

\author[1]{Quanjiang Yu\corref{mycorrespondingauthor}} 
\address[1]{Department of Mathematical Sciences, Chalmers University of Technology and  University of Gothenburg, SE-42196 Gothenburg, Sweden}
\ead{yuqu@chalmers.se}
\author[2]{Pramod Bangalore}
\address[2]{Greenbyte AB, SE-411 09 Gothenburg, Sweden} 
\author[3]{Sara Fogelstr{\"o}m}
\address[3]{Department of Electrical Engineering, Chalmers University of Technology } 
\author[1]{Serik Sagitov}



\cortext[mycorrespondingauthor]{Corresponding author}


\begin{abstract}
We suggest a mathematical model for computing and regularly updating the next preventive maintenance plan for a wind farm. Our optimization criterium takes into account the current ages of the key components, the major maintenance costs including eventual energy production losses as well as the available data monitoring  the condition of the wind turbines. We illustrate our approach with a case study based on data  collected from several wind farms located in Sweden. Our results show that preventive maintenance planning gives some effect, if the wind turbine components in question live significantly shorter than the turbine itself.
\end{abstract}

\begin{keyword}
 Preventive maintenance \sep Linear programming \sep Cox proportional hazards \sep Wind turbine \sep Weibull survival function
\end{keyword}

\end{frontmatter}


\section{Introduction}

Renewable energy sources like Wind and Solar are set to play a major role in the energy systems of the future. According to some projections, like the one presented in \cite{energyoutlook}, more than 50\% of total electricity might come from renewable energy sources by 2050. These projections mean that the wind turbines will grow to a much larger number in the near future, both on and offshore. To cope with such large number of assets, it will become crucial to automate processes around operation and maintenance of these systems. In addition to simplifying and streamlining decision making, such automated processes might also allow for optimisation around maintenance costs, which even today account for quite a considerable portion of the operational life cycle cost for wind turbine assets, especially for offshore installations. 

Maintenance optimisation can be viewed as the process of deciding the best moment in time, both from economic and technical perspectives, to replace one or more components considering the impact of each maintenance activity on the life cycle cost of an asset or assets. The recent literature on wind turbine preventive maintenance planning extends the modelling scope by paying special attention to particular performance factors for the wind power systems. Paper \cite{zheng2020optimal} looks into the effects of the varying wind speed on the wind turbine maintenance planning. Paper \cite{davoodi2020preventive} singles out the converter as a crucial component of the wind turbine and builds an optimization model to find the optimal replacement times for the converters. Papers \cite{wang2020optimal} and \cite{zhang2017opportunistic} deal with imperfect preventive maintenance. 

By utilizing the information about the state of various critical components the maintenance routines can be further improved. Cox's Proportional Hazards Model (PHM), proposed in \cite{cox1984analysis}, utilizes measurable entities as covariates to update the hazard function for a component, making the PHM quite handy for application with data from a condition monitoring system (CMS). Several research teams have suggested various optimisation models in an attempt to make use of condition monitoring data by applying some version of the PHM, see for example, 
\cite{wu2011optimal}, \cite{you2011updated}, \cite{li2013multi}, \cite{pham2012machine}. Furthermore, Paper \cite{vlok2002optimal} developed a probabilistic model to estimate the remaining lifetime of a machinery using data from a CMS. Their probabilistic approach involves a PHM with Weibull  baseline hazard and a Markov process model. Vibration data is used as an input from the CMS to illustrate a practical application of this probabilistic model. 
Similarly in \cite{banjevic2006calculation}, the failure process along with the covariate process is represented by a discrete Markov process. A PHM algorithm is proposed for predicting the remaining lifetime of the machinery based on a condition monitoring process.

In \cite{wang2015reliability}, the authors feed the online vibration and temperature signals of bearings from the CMS into a neural network and predict the features of bearing vibration signals at any time horizon. Furthermore, according to the features, degradation factor was
defined. A PHM was generated to estimate the survival function and forecast the remaining lifetime of the bearing.

Paper \cite{ghasemi2009evaluating}  is built upon a hidden Markov model, assuming that the
equipment’s unobservable degradation state evolves as a
Markov chain. The
Bayes rule is used to determine the probability of being in a certain
degradation state at each observation moment. Cox’s time-dependent PHM is applied to deal with the equipment’s failure rate. Two main problems are addressed: the problem of imperfect observations, and the problem of taking into account the whole history of observations.

The recent papers \cite{bangalore2018analysis}, \cite{bangalore2017artificial}  develop a machine learning approach to maintenance scheduling for a wind turbine whose condition is monitored by a time series 
$\{x(1), x(2),\ldots, x(t)\}$
summarising some key characteristics of the turbine which can be used for predicting the failure times after time $t$. A deep learning algorithm was trained for a prediction $\hat x(t+1)$ of the next value $x(t+1)$ based on a time series observed up to the current time $t$. Then at time $t+1$, depending on  a certain measure of discrepancy between the observed $x(t+1)$ and predicted $\hat x(t+1)$ values,
a decision is made whether a PM should be performed in the near future or not.
One of the key simplifying assumptions requires that the turbine's component in question has an exponential life length distribution.

Most of the research towards condition based maintenance utilizes the data from vibration based CMS. The vibration data is measured at a high frequency, in range of \textit{kHz}, and the data is processed using various algorithms before it is stored. It might, in some cases, also be difficult to get access to data from the vibration based systems. Hence, in this paper we focus on creating a decision making model which utilizes easily available signals such as component temperatures. 

Modern wind turbines are equipped with Supervisory Control and Data Acquisition (SCADA) system which measures and stores the data for various component temperatures, this data has been utilized along with information about historical failures to create a model to estimate condition based failure rate of gearboxes. Furthermore, in this paper, the common assumption of exponential life length distribution (constant failure rate) is relaxed using the Weibull statistical model. 

The optimisation model presented in this paper is a slight variation of the one proposed in \cite{Paper}, where a multiple component setting for a single wind turbine (without  condition monitoring data) was addressed.

The rest of the paper is organized as follows. Section \ref{CW} describes how utilising condition monitoring data one can estimate the Weibull parameters of aging gearboxes. Section \ref{modmodel} gives a detailed description of the optimisation model for $n$ wind turbines each represented by its gearbox. 
Section \ref{case} presents a detailed case study based on data collected from several onshore wind farms in Sweden. Section \ref{test} has a closer look at a particular wind farm. Finally, Section \ref{conc} contains discussions and conclusions.

\section{Weibull parameters under condition monitoring}\label{CW}

The key ingredient of the optimisation model of this paper is the Weibull distribution for the life length $L$ of a generic gearbox
\[\rP(L>t)=e^{-\theta t^{\kappa}},\quad t>0.\]
It is assumed that under the normal conditions the Weibull parameters $(\theta,\kappa)$ of a gearbox take certain baseline values $(\theta_0,\kappa_0)$, so that the baseline hazard function (failure rate at age $t$) takes the form
\[r_0(t)=\theta_0 \kappa_0 t^{\kappa_0-1}.\]
Suppose that time series data 
$$\boldsymbol x=(x(1), x(2),\ldots)$$
measuring an appropriate covariate $x(t)$ at different times $t$, can be utilized to monitor the condition of a given gearbox. 
Assuming that the shape parameter $\kappa_0$ of the Weibull distribution of the gearbox's life length is constant over time,  
the task addressed in this section is to find an appropriate value of the scale parameter $\theta=\theta(\boldsymbol x,t)$ which would update the failure rate of the gearbox 
\begin{equation}
r(\boldsymbol x,t)=\theta(\boldsymbol x,t) \kappa_0 t^{\kappa_0-1}
\label{rex}
\end{equation}
by incorporating the available condition monitoring data $x(1),\ldots,x(t)$. 

\subsection{Finding $(\theta_0, \kappa_0)$ using training data}\label{trda}

Consider a set of historical data containing the observed ages of  still operational gearboxes $u_1,\ldots,u_K$, and yet another historical data set for gearboxes that have failed
\begin{align}
(v_1,\boldsymbol x^{(1)}),\ldots,(v_N,\boldsymbol x^{(N)}),
 \label{dest}
\end{align}
where, $v_k$ is the  failure age of a gearbox $k$, and $\boldsymbol x^{(k)}=(x_1^{(k)},\ldots,x_{v_k}^{(k)})$ is the corresponding recorded history of the monitoring data.
The baseline Weibull parameter values $(\theta_0,\kappa_0)$ are estimated from the two sets of observed life times
\[U=\{u_1,\ldots,u_K\},\quad V=\{v_1,\ldots,v_N\},\]
by maximising the likelihood function
\begin{equation*}
\mathcal L(\theta,\kappa)=\prod_{t\in V}\rP(L=t)\prod_{t\in U}\rP(L>t)=\prod_{t\in V}(e^{-\theta (t-1)^\kappa}-e^{-\theta t^\kappa})\prod_{t\in U}e^{-\theta t^\kappa}.
\end{equation*}

\subsection{The use of proportional hazard method }\label{trda}

To be able to update the hazard rate by means of \eqref{rex}, the following version of PHM is suggested:
\begin{align}
 \theta(\boldsymbol x,t)=\theta_0 e^{\beta(\bar x(t)-\bar x)},\quad j=1,\ldots,n,
 \label{est}
\end{align}
where
\[\bar x=\frac{x(1)+\ldots+x(12)}{12}\]
is the first year average of the covariate $x$
and 
\[\bar x(t)=\frac{x(t-2)+x(t-1)+x(t)}{3}\]
is the latest three-month moving average.
Obviously, this approach requires that the farm has been in operation for at least $15$ months.

The Cox regression parameter $\beta$ mentioned in \eqref{est} is estimated from the training data set \eqref{dest} assuming that the data is labeled in such a way that the failure times are sorted in the ascending order
\[v_1<v_2<\ldots<v_N.\]
The key argument of the Cox method  \cite{cox1984analysis} is that expressions \eqref{rex} and \eqref{est} imply the following expression for the partial likelihood function of the regression parameter $\beta$
\[\mathcal L^*(\beta)=\prod_{j=1}^{N}\frac{r(v_j,\boldsymbol x^{(j)})}{\sum_{i=j}^{N} r(v_j,\boldsymbol x^{(i)})}=\prod_{j=1}^{N}\frac{\exp\{\beta\bar x^{(j)}(v_j)\}}{\sum_{i=j}^{N} \exp\{\beta\bar x^{(i)}(v_j)\}}.\]
Maximisation of the partial likelihood $\mathcal L^*(\beta)$ leads to the desired maximum likelihood estimate $\beta_0$.

As a result for the current $n$ component setting, we obtain the updating formulas for the $n$ pairs of the Weibull parameters
\begin{align}
\theta_j=\theta_0\phi_j(t), \quad \kappa_j=\kappa_0, \quad j=1,\ldots,n,
\label{thest}
\end{align}
involving Cox factors 
\begin{align}
\phi_j(t)=e^{\beta_0(\bar x_j(t)-\bar x_j)}, \quad j=1,\ldots,n,
 \label{west}
\end{align}
based on $n$ times series
\[\boldsymbol x_j=(x_j(1), x_j(2),\ldots),\quad j=1,\ldots,n.\]

The Cox factor \eqref{west} has the following effect on the failure rate of the gearbox $j$, provided $\beta_0$ is positive (in other words, assuming that the chosen covariate is such that higher values of $x(t)$ indicate  higher stress on the gearbox at time $t$). At the time of observation $t$, the first year average $\bar x_j$ is compared with the last three month average $\bar x_j(t)$. 
If the difference $\bar x_j(t)-\bar x_j$ is close to zero, then the current condition of the turbine $j$ is deemed to be normal and 
formula  \eqref{thest} suggests using the baseline parameters $\theta_j=\theta_0, \kappa_j=\kappa_0, $ for describing the failure rate of the gearbox $j$. 
However, if it turns out that $\bar x_j(t)>\bar x_j$, so that $\theta_j>\theta_0$, then the corresponding hazard rate
\[r_j(t)=\theta_j \kappa_0 t^{\kappa_0-1}\]
becomes larger that the base line value $r_0(t)$.
Alternatively, if $\bar x_j(t)<\bar x_j$, then of course, the failure rate of the gearbox at time $t$ is below the normal: $r_j(t)<r_0(t)$.

\section{Optimal Preventive Maintenance schedule for $n$ gearboxes}\label{modmodel}

An efficient optimisation model for a single wind turbine with several components was presented in \cite{Paper}. In this section the optimisation model from \cite{Paper} is adapted to a setting with $n$ wind turbines, where each wind turbine is represented by its gearbox as the key component. Section \ref{smodel} introduces the main cost parameters including so-called virtual maintenance costs. Section \ref{npm} presents the main step of our optimal scheduling algorithm summarized in Section \ref{reop}.

\subsection{Maintenance costs}\label{smodel}
The maintenance costs of gearboxes are modelled in terms of the following parameters
\begin{description}
\item[ $g$] is the total cost of a corrective maintenance (CM), including the logistic cost, down-time cost, and the cost of a new gearbox;
\item[ $h_0$] is the fixed cost of a preventive maintenance (PM) activity, this cost is the same regardless of how many gearboxes are planned to be replaced during this activity;
\item [$h$] is the variable cost related to the PM replacement which takes into account the replacement cost of a gearbox, the downtime cost, and the initial  value loss of the gearbox in use;
\item [$m$] is the monthly loss of the value for a gearbox in use. 
\end{description}
To illustrate the use of the parameters $(h_0,h,m)$, consider a PM plan suggesting to simultaneously replace three components having ages $(a_1,a_2,a_3)$ in months. Then the total cost associated with this PM activity, $f$, is calculated as
\begin{equation*}
 f\ =\ h_0+(h+a_1m)+(h+a_2m)+(h+a_3m)\ = \ h_0+3h+(a_1+a_2+a_3)m. 
\label{eq:pmcost}
\end{equation*}

Given the Weibull parameter values $(\theta_j,\kappa_0)$, using the approach of \cite{Paper}, the virtual replacement cost $b_j(a)$ for the gearbox of  age $a$ can be computed.
(For further details on the exact calculation of $b_j(a)$ and the interpretation of the virtual replacement cost based on the renewal-reward argument, the reader is referred to \cite{Paper}.)
In what follows, 
\[B_j(a)=(h+am)\wedge b_j(a)\]
stands for the minimum between two age specific costs: the age-specific PM cost and the virtual replacement cost. 

\subsection{Monthly maintenance replacement cost $c$}\label{Lgm}

Consider a wind farm with $n$ new gearboxes at time $t=0$, where time to failure of the first gearbox is denoted by
\[L_0=\min(L_1,\ldots, L_n)\]
By independence, we have 
\[\rP(L_0>t)=\prod_{i=1}^n\rP(L_i>t)=e^{-n\theta t^{\kappa}}.\]
If the next PM is planned at time $t$, then the first renewal time of the system $X=X(t)$ can be calculated as presented in Equation \ref{Xt}.
\begin{equation}
X=L_0\wedge t=L_0\cdot 1_{\{L_0\leq t\}}+ t\cdot 1_{\{L_0> t\}} 
\label{Xt}
\end{equation}
The corresponding reward value $R=R(t)$, can be computed as 
%
\[R=(g+(n-1) B_0(L_0)) 1_{\{L_0\leq t\}}+ (h_0+nB_0(t))1_{\{L_0> t\}},\]
where 
\[B_0(a)=(h+am)\wedge b_0(a)\]
is the age specific replacement cost, provided the gearbox's Weibull parameters take the baseline values. 

Then, the  renewal-reward theorem implies that the time-average maintenance cost $\frac{\rE(R)}{\rE(X)}$ is  the following function of the planning time $t$:
\[q_t=\frac{g\rP(L_0\le t)+(n-1)\rE( B_0(L_0)\cdot 1_{\{L_0\leq t\}})+(h_0+nB_0(t)) \rP(L_0> t)}
{\rE(L_0 \cdot 1_{\{L_0\leq t\}})+ t\rP(L_0> t)}.\]
After minimising $q_t$ over $t$ we can define the monthly maintenance replacement cost of the wind farm as a constant
\begin{equation}
c=\min_{t\ge 1}q_t \label{c}, 
\end{equation}
see \cite{Paper} for a more detailed explanation.

\subsection{The key optimization step}\label{npm}

For the planning period $[s,T]$, where $T$ is the end of life for the whole wind farm,a PM plan any array can be defined as
$$(\boldsymbol w_s,\boldsymbol y_s,z)=\{w^{j}_{t},y_t,z: \ 1\le  j\le n,\ s+1\le t\le T\}$$  
with binary components $w^ {j}_{t},y_{t},z\in \{0,1\}$ satisfying the following linear constraints
\begin{subequations}
\begin{align}
y_t&\geq w^{j}_{t},\quad t=s+1,\ldots, T,\ j=1,\ldots,n,\label{tu}\\
    \sum_{j=1}^n w_t^j&\geq y_t,\quad t=s+1,\ldots,T,\label{x_t^j}\\
\sum_{t=s+1}^{T} y_t&=1-z.\label{y}
\end{align} 
\label{const}
\end{subequations}
Here, $w^ {j}_{t}=1$ means that at time $t$  a PM activity is planned for turbine $j$, otherwise $w^ {j}_{t}=0$. 
Similarly,
 $y_{t}=1$  means that at time $t$ a PM activity is planned for at least one of the turbines in the wind farm, constraints \eqref{tu} and \eqref{x_t^j}. The equality $z=1$ means that no PM activity is planned during the whole time period $[s+1,T]$, constraint \eqref{y}.

Given the ages of $n$ gearboxes at time $s$ 
\[\boldsymbol a=(a_1,\ldots,a_n),\]
the first failure time is $s+L_{\boldsymbol a}$, where (lifting the turbine index upstairs)
$$L_{\boldsymbol a}=\min(L^1_{a_1},\ldots, L^n_{a_n}).$$
and
\[\rP(L_a^j> t)=\exp\left\{\theta_j\big(a^{\kappa_0}-(a+t)^{\kappa_0}\big)\right\},\quad t\ge0,\]
is the survival function conditional on the age $a$.

The cost assigned to a PM plan can be denoted as
\begin{align*}
    F_{(s,\boldsymbol a)}(\boldsymbol y_s,z)
   & =\sum_{t=s+1}^{T}\Big( g+(T-s-L_{\boldsymbol a}) c+\sum_{j\neq \gamma} B_j(a_j+L_{\boldsymbol a})\Big)1_{\{s+L_{\boldsymbol a}\leq t\}}y_t \\
   &+\sum_{t=s+1}^{T} \Big(  h_0+(T-t)c+\sum_{j=1}^{n}B_j(a_j+t-s)\Big)1_{\{s+L_{\boldsymbol a}>t\}}y_t\\
   &+\Big( g+(T-s-L_{\boldsymbol a}) c+\sum_{j\neq \gamma} B_j(a_j+L_{\boldsymbol a})\Big)1_{\{s+L_{\boldsymbol a}\leq T\}}z,
\end{align*}
where $\gamma$ is the label of the gearbox that failed at time $s+L_{\boldsymbol a}$.
Notice that the total cost function $F_{(s,\boldsymbol a)}(\boldsymbol y_s,z)$ does not explicitly depend on $\boldsymbol w_s$. The role of $\boldsymbol w_s$ becomes explicit through the following additional constraint 
\begin{align}
    &(h+(a_j+t-s)m)\cdot w_t^j+b_{j}(a_j+t-s)\cdot (y_t-w_t^j)=B_j(a_j+t-s)\cdot y_t, \label{x_j}\\
    & \quad \quad t=s+1,\ldots, T,\quad \ j=1,\ldots,n.\nonumber
\end{align}
If  $y_t=1$, that is if a PM activity for at least one component is scheduled at time $t$, then for each component $j$, there is a choice between two actions at time $t$: either perform a PM, so that $w_t^j=1$ and $y_t-w_t^j=0$, or do not perform a PM and compensate for the future extra costs caused by the current gearbox age using the virtual replacement cost value (corresponds to $w_t^j=0$ and $y_t-w_t^j=1$).

The optimal maintenance plan according the presented approach is the solution of the linear optimisation problem
\begin{align*}
 \text{minimise}\qquad & f_{(s,\boldsymbol a)}( \boldsymbol y_s,z)=\rE(F_{(s,\boldsymbol a)}(\boldsymbol y_s,z))\\
\text{subject to}\qquad & \text{linear constraints } \eqref{tu},\eqref{x_t^j}, \eqref{y}, \text{ and }\eqref{x_j},\\
&w_t^j\in \{0,1\},\quad t=s+1,\ldots T,\ j=1,\ldots,n,\\
&y_t\in \{0,1\},\quad t=s+1,\ldots T,\\
&z\ \in\{0,1\}.
\end{align*}

\subsection{Optimal scheduling algorithm for $n$ gearboxes}\label{reop}

\begin{algorithm}
\caption{Optimal scheduling algorithm}
{
\small
 Input:  $s,\ T,\ g,\ h_0,\ h,\ m,\ \kappa_0,\ \beta_0,\ \theta_0,\ a_1,\ldots,a_n$
 
Step 1: \textbf{for} $j=1:n$

\hspace{1.65cm}\textbf{if} {$a_j\leq2$} \textbf{then}

\hspace{2.0cm}Set $\theta_j:=\theta_0$

\hspace{1.65cm}\textbf{else}

\hspace{2cm}Collect the last three months of condition monitoring data 

\hspace{2cm}and compute $\theta_j$ based on $(x_j(s-2),\ x_j(s-1),\ x_j(s))$

\hspace{1.65cm}\textbf{end if}

\hspace{1.15cm}\textbf{end for}

Step 2: Apply the key optimization step, see Section \ref{npm}, with 

\hspace{1.15cm}Output $t^*,\ \mathcal{P}\subset\{1,\ldots,n\}$


Step 3:  Suppose after time $s$, the first failure would be at time $t'$

\hspace{1.15cm}\textbf{if} $t' \leq\min \{t^*, s+3\}$ \textbf{then}

\hspace{1.5cm}Put $t^*:=t'$

\hspace{1.5cm}Go to Step 4




\hspace{1.15cm}\textbf{else}

\hspace{1.5cm}Go to Step 5

\hspace{1.15cm}\textbf{end if}

Step 4: Apply opportunistic maintenance step at time $t^*$ with 

\hspace{1.15cm}Output $\mathcal{P}\subset\{1,\ldots,n\}$


\hspace{1.15cm}Go to Step 6

Step 5: \textbf{if} $t^*\leq s+3$



\hspace{1.5cm}Go to Step 6

\hspace{1.15cm}\textbf{else}

\hspace{1.5cm}Update $a_j:=a_j+3,\ j\in\{1,\ldots,n\}$; $s:=s+3$

\hspace{1.5cm}Go back to Step 1

\hspace{1.15cm}\textbf{end if}

Step 6: The gearboxes with labels in $\mathcal{P}$ are replaced by new ones

\hspace{1.15cm}Update $a_j:=0,\ j\in\mathcal{P}$; $a_j:=a_j+t^*-s,\ j\notin\mathcal{P}$; $s:=t^*$

\hspace{1.15cm}Go back to Step 1

}
\label{al1}
\end{algorithm}

In this section, the main result of this paper is summarized in the form of Algorithm \ref{al1} producing a PM plan for a given planning period $[s,T]$,
focusing on the gearbox components of $n$ wind turbines constituting a wind farm. It is assumed that the starting planning time $s$ is such that $s\ge15$, and that the length of the updating period is $3$ months. The following data and parameters are assumed to be available:
\begin{description}
\item[ --] condition monitoring time series $\boldsymbol x_j=(x_j(1), x_j(2))\ldots$ for $j=1,\ldots,n$,
\item[ --] baseline Weibull parameters $(\theta_0,\kappa_0)$ and Cox regression parameter $\beta_0$ obtained from the training data,
\item[ --] maintenance cost parameters $(g,h_0,h,m)$,
\item[ --] gearbox ages $\boldsymbol a=(a_1,\ldots,a_n)$ at time $s$.
\end{description}

The key step of Algorithm \ref{al1}, Step 2, is described in Section \ref{npm}.

Step 4 requires clarification. If any of the gearboxes breaks down before the planned next PM, a CM replacement alongside with opportunistic replacements are performed. The opportunistic replacement work as follows: since the maintenance personal need to go there and perform CM on the broken component, they may as well maintain other gearboxes if they are close to break down to save the logistic cost. So, for each other component, we compare the virtual maintenance cost and the PM cost, if virtual maintenance cost is higher, it means that the gearbox is too old, it is beneficial to perform PM on the corresponding gearbox, more details see section 7 in \cite{Paper}. 
After each replacement (either PM or CM) one has to update the vector of ages and the starting time $s$ accordingly, and then repeat the key step of the algorithm.

\begin{center}
\begin{table}[ht!]
  \begin{tabular}{l|c|c|c|c}
 Farm ID &\# turbines&\# failures&Time in use (month)&SCADA   \\\hline
  1& 5&0&126&yes\\
  2& 5&2&143&yes\\
  3&6&1&124&yes\\
  4& 8&8&146&yes\\
  5& 8&10&140&yes\\
  6& 9&0&72&yes\\
  7&11&0&72&yes\\
  8&13&0&72&yes\\
  9&16&8&137&yes\\
  10&1&0&101&no\\
  11&3&0&113&no\\
  12&5&0&114&no\\
  13&5&5&168&no\\
  14&6&0&150&no\\
  15&9&0&94&no\\
  16&10&1&79&no\\
  17&10&0&115&no\\
  18&12&7&144&no\\
  19&12&0&91&no\\
  20&32&4&95&no
  \end{tabular}
   \caption{Summary of the data on 20 farms: the 2nd column gives the number of wind turbines constituting each farm, the 3rd column gives the number of gearbox failures during the number of months mentioned  in the 4th column. The 5th column says whether a farm has temperature sensor data in SCADA or not.}
 \label{tabu1}
  \end{table}
\end{center}

\section{Swedish data set on 20 farms}\label{case}

The case study is based on data collected in November 2020 on 20 wind power farms located in Sweden, see Table \ref{tabu1}. The wind farms are located in southern and middle part of Sweden and were erected from 2006 to 2014 (1 from 2006, 3 from 2008, 3 from 2009, 2 from 2010, 3 from 2011, 1 from 2012, 3 from 2013, and 4 from 2014) . Column 1 sets labels to the farms, column 2 gives the number of turbines in each farm,  column 3 gives the observed number of gearbox failures for the respective farm, column 4 says during how many months the farm was observed, and finally, column 5 specifies whether the farm has temperature sensor data in SCADA or not. For example, the data for wind farm 9 that has 16 turbines is available for 137 months, and during this period of time, the wind farm experienced 8 gearbox failures. A detailed case study on the data from wind farm 9 is presented in Section \ref{test}.

The total number of turbines is 186. The total number of gearboxes, in this data set is 232, with 46 gearboxes that have failed and 186 gearboxes that still are in use. There are $5$ wind turbines that has experience gearbox break down twice. Using the method described in Section \ref{trda}, one arrives to the following baseline parameter values
\[\theta_0=8.386\cdot 10^{-4}, \quad \kappa_0=1.217,\]
corresponding to the mean life length for a gearbox of $316$ months or 26 years. This estimate is in contrast to the reliability analysis results presented in the literature reporting much shorter life lengths for the gearboxes. However, the result is not surprising given that the data set consists only of onshore and relatively new wind turbines. Over the years lot of progress has been made in design of wind turbine gearboxes which has lead to fewer failures in more stable conditions. Furthermore, it must be noted that certain wind farms in the study have had an unusually high number of gearbox failures; for example wind farms 5 and 9. The estimated life expectancy for gearboxes in these wind farms is much shorter than $316$ months. The maintenance optimization method presented here is beneficial when the life expectancy of the gearboxes is much shorter than the planned life of the wind turbines. Hence, in order to demonstrate the applicability of the method, the case studies are based on the following Weibull parameter values, presented in \cite{tian2011condition},
\begin{equation}
\theta_0=1.95\cdot 10^{-6}, \quad \kappa_0=3.\
\label{tuk}
\end{equation}
With these baseline Weibull parameters, the mean life length for a gearbox becomes $71$ months.

\begin{center}
\begin{table}[ht!]
  \begin{tabular}{l|c|c|c|c}
$k$&Farm ID& Failure time $a_k$&Cox factor $\phi_k$&$\phi(a_k)$    \\
 \hline
1&4&  21&2.60    & $0.97\pm0.12$ \\
2&4&  25& 2.72   & $1.16\pm0.11$  \\
3&4&  25& 3.47   & $1.16\pm0.11$ \\
4&9&  25& 3.61   & $1.16\pm0.11$ \\
5&4&  37& 3.85   & $0.96\pm0.14$ \\
6&5&  37&  1.40  & $0.96\pm0.14$ \\
7&9&  43&  1.67  & $0.98\pm0.10$ \\
8&2&  52&  1.20  & $1.15\pm0.08$ \\
9&5&  61&  2.12  & $0.91\pm0.11$ \\
10&5&  61&  2.06  & $0.91\pm0.11$ \\
11&5&  66&  1.26  & $1.25\pm0.16$ \\
12&4&  73&  0.91  & $1.07\pm0.25$ \\
13&9&  73&  1.42  & $1.07\pm0.25$ \\
14&9&  73&  1.33  & $1.07\pm0.25$ \\
15&5&  80& 1.80   & $0.99\pm0.19$ \\
16&2&  97&2.04& $1.01\pm0.29$\\
17&4&  97&0.97& $1.01\pm0.29$\\
18&9&  97&1.04& $1.01\pm0.29$\\
19&9&  109&1.28&$1.09\pm0.30$\\
20&3&  116&1.56&$1.03\pm0.32$\\
21&9&  121&1.71&$1.06\pm0.31$\\
22&9&  121&1.23&$1.06\pm0.31$\\
23&5&  133&0.93&$1.10\pm0.39$
\end{tabular}
  \caption{The data on 23 SCADA connected gearboxes that went down during the time of observation. Column 2 specifies at which farm the failure was observed. Column 3 gives the life length of the gearbox. Column 4 gives gearbox specific Cox factor at the time of failure. Column 5 gives the Cox factors averaged across 55 non-failed gearboxes at the matching ages.
}
 \label{tab2}
  \end{table}
\end{center}

According to Table \ref{tabu1}, among the gearboxes for which the SCADA condition monitoring data is available, 29 have experienced a failure. Out of these 29  gearboxes, 4 belonged to wind farm 4, which has been connected to SCADA since month 52 of exploitation time, furthermore, 2 of the failure times did not satisfy the requirement of 15 months monitoring data available. This leaves us with 23 gearboxes to which our approach can be applied. 

Table \ref{tab2} focuses on  $23$ gearboxes whose failure times are given in column 3 and for which the SCADA monitoring data is available. Implementing the approach of Section \ref{trda} based on \eqref{tuk}, we applied the steepest descent algorithm and obtained
\begin{equation}
\beta_0=0.203.
\label{bee}
\end{equation}
 Column 4 of the Table \ref{tab2} gives the Cox factors $\phi_k=\phi_k(a_k)$ obtained using \eqref{west} with $j=k$ and $t=a_k$, that is at the time prior to the failure of the gearbox in question. An immediate observation is that $19$ out of $22$ values $\phi_k$ are higher than the critical value 1, an indication of the increased risk of failure (conditioned on the current age). However, these results are very sensitive to the estimate $\beta_0$. It is more relevant to compare the Cox factor of the failed gearbox to the gearboxes which were still functioning at the  age given in column 3, see column 5 containing 95\% confidence intervals. For the majority of gearboxes in use, the Cox factor  $\phi(a_k)<\phi_k$ is estimated to be smaller.

\section{Wind farm 9: test study using historical data}\label{test}
Here we use the historical data available for wind farm 9 to see if our approach, based on estimates \eqref{tuk} and  \eqref{bee}, is able to avoid the failure events by placing PM activities at right times and for the right gearboxes.
Recall that wind farm 9 consists of $16$ wind turbines, with 8 of them having experienced failures at ages given in the table below.
\begin{center}
\begin{tabular}{l|cccccccc}
Gearbox ID index	$k$	&4	&7 	&13	&14	&18	&19	&21	&22\\\hline
Failure time (months)	&25	&43	&73		&73	&97	&109&121		&121

\end{tabular}
\end{center}
Observe that two pairs of equal failure times indicate violations of the model assumption of independence between the gearbox life times. Our guess is that for each of the paired events,  one of the gearboxes might have broken down earlier and the turbine stayed idle until the second gearbox went down, so that both gearboxes were replaced simultaneously.

\begin{figure}[h]
\centering
\includegraphics[width=0.8\textwidth]{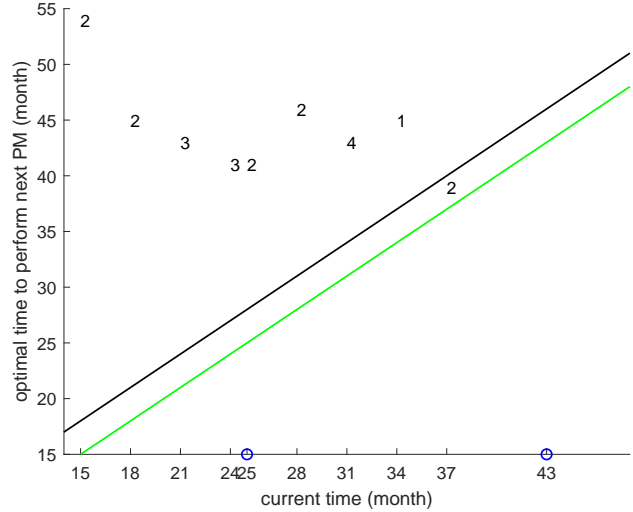}

        \caption{The recurrent next PM planning for wind farm 9}
            \label{cnm}
\end{figure}

The results of our study based on the historical data for wind farm $9$ are summarised in Figure \ref{cnm}. It shows the recurrent $3$-month updates of the PM planning, so that if the next PM activity is planned later than in the next $3$ month time period, it will not be performed. After $3$ months, we update the data from the CMS and resolve the optimal problem again to obtain a new maintenance plan. The green line represents the observation time and the black line represents the planning horizon three months ahead. Each planning round giving the next time for PM as a point lying above the black diagonal, will be followed by a new planning round with an updated time for the  next PM. The next PM plan will be implemented only if the next PM point lies between the two diagonals on the graph.

As shown on  the $x$-coordinate, the first PM schedule was produced at time $15$. 
The resulting optimal planning time at month $54$ is shown on the $y$-coordinate. The corresponding point $(15,54)$ is marked on the graph by label $2$ telling that $2$ gearboxes out of $16$ should be replaced at time $54$. Since point $(15,54)$ lies above the black diagonal, we apply our algorithm once again at time $15+3=18$ and find the new PM time to be at month $45$ when $2$ gearboxes should be replaced. At time $21$, an updated PM plan says that $3$ gearboxes should be replaced at month $43$.  And so on.

The most interesting points on the graph are times $25$ and $37$. 
At time $24$, the optimal PM plan was to replace gearbox $4$ at month $41$ together with two other gearboxes. What happened next, according to the historical data, is that gearbox $4$ broke down at time $25$. Since we need to perform CM on gearbox $4$, we apply opportunistic maintenance. We found out for all other gearboxes, the virtual maintenance cost for each gearbox is  lower than the corresponding PM cost. Thus, the optimal plan at month $25$ before the replacement is to only perform CM at gearbox $4$.  After the CM, we resolve the optimal model with update data of gearbox $4$, i.e., age $0$ and baseline Weibull parameters. The optimal plan is to replace $2$ gearboxes at month $41$.

First at time $37$, the next PM time falls within the three month window. It means that in a planned manner $2$ gearboxes should be replaced at month $39$. From the historical data, we see that gearbox $7$ that has failure time $43$ in the data, is closest to this PM time and is among the PM replaced gearboxes.

In this case study, we used the following parameter values for the various maintenance costs, we normalized the data and use a virtual monetary unit:
\begin{equation}
   g=1+d_t, \quad m=0.008, \quad h_0=0.13,\quad h=0.294+d_t/6, 
\label{tre}
\end{equation}
where the downtime cost $d_t$ depends on the month of the replacement:
\begin{center} 
  
  \begin{tabular}{l|cccccccccccc}
Month&Jan&Feb&Mar&Apr&May&Jun
\\\hline 
$d_t$&0.075&0.044&0.067&0.053&0.059&0.069
\end{tabular}

  \begin{tabular}{l|cccccccccccc}
Month
&Jul&Aug&Sep&Oct&Nov&Dec\\\hline 
$d_t$
&0.046&0.070&0.085&0.066&0.066&0.057
\end{tabular}
\end{center} 
The monthly downtime cost is calculated from monthly productions multiplied with monthly selling price and which averaged over three years. The production for each month comes from data from the eight turbines in wind farm 9 that haven't replaced gearboxes yet and is from year 2017 to 2019. The month selling price is a combination of monthly electricity spot prices from Nord Pool and monthly prices for the green certificates from Svensk kraftmäkling for the same three years.
 
The CM cost is
\begin{equation}
    g=c_g+c_m+d_t
\end{equation}
where $c_g$ is the cost for a new gearbox (0.64) and $c_m$ is the maintenance cost. The maintenance cost is divided into four parts: transport cost for the crane (0.04), set-up cost for the crane (0.09), working cost for the crane (0.16) and the manpower cost for replacing a gearbox (0.07). The sum of $c_g$ and $c_m$ is 1 virtual monetary unit. 
The different costs comes from three different wind power operators and the presented costs are averaged and normalized values from their data. 

The shared maintenance cost for PM, $h_0$, consists of two parts of the maintenance cost; the transport cost for the crane (0.04) and the set-up cost for the crane (0.09), bringing it to a total of 0.13 virtual monetary unit. 

According to specification \eqref{tre}, $h$ consists of the other two parts of the maintenance cost; the working cost for the crane (0.16) and the manpower cost for replacing a gearbox (0.07), as well as the initial loss of value of the gearbox (0.064) and the downtime cost for PM. The total is 0.294 plus downtime cost. The initial loss of value of the gearbox is set to 10\%  of the value of a new gearbox.

The monthly value loss $m$ is set to approximate 0.008 and is defined as the gearbox cost (0.64) minus the initial value loss (0.064) divided by the expected life time of the gear box (71 months \cite{tian2011condition}) in virtual monetary unit:
\begin{equation}
 m=\frac{0.64-0.064}{71}\approx0.008. 
\end{equation}

We assume an equal depreciation of the value of the gearbox per year during its lifetime. 

Notice that the PM down time cost is $6$ times smaller than that of the CM counterpart, since a PM activity goes $6$ times faster. 


\section{Conclusions}\label{conc}
In this paper, we adapted the optimisation model of \cite{Paper} developed for a single wind turbine with $n$ components to a setting with $n$ wind turbines constituting a wind farm. Then, the model was enhanced by adding a parameter updating step allowing our maintenance scheduling optimisation algorithm to take into account the real time data from the CMS. This parameter updating  step is based on the Cox proportional hazards method.

Using the suggested approach, we studied the recent historical data from several wind farms located in Sweden.  A more careful analysis was performed using the data from one of these farms. Our analysis showed that the success of the scheduling using our model to high extend depends on the baseline values of the Weibull parameters. One of the clear conclusions of our analysis is that PM planning gives some effect, only if the wind turbine components in question live significantly shorter than the turbine itself.
Provided the component's life time is notably shorter than the turbine's life time, our approach may result in appreciable savings due to smart scheduling of PM activities by monitoring the ages of the components in use as well as using available real time data supervising the condition of the wind turbines in a wind farm.

\section*{Acknowledgements}
\noindent We acknowledge
the financial support from the Swedish Wind Power Technology Centre.


\bibliography{mybibfile}

\end{document}